\input amstex
\documentstyle{amsppt}
\magnification=1200
\hcorrection{0.5in}
\hsize =14truecm
\baselineskip =18truept
\vcorrection{0.15in}
\vsize =20.5truecm

\nologo

\NoRunningHeads

\NoBlackBoxes

\topmatter
\document

\title {\bf Rational connectedness of log $Q$-Fano varieties}
\endtitle
\author  {Qi  Zhang}    
\endauthor
\address{Department of mathematics,  University of
Missouri, Columbia, 
MO
65211. U.S.A}  \endaddress
\email{qi\@math.missouri.edu}\endemail
\abstract{In this paper, we give an affirmative answer to a conjecture
in the Minimal Model Program. We prove that log $Q$-Fano varieties are
rationally  connected. We also study the behavior of the canonical
bundles under projective morphisms}
\endabstract
\endtopmatter

\heading{\bf $\S$1.    Log  $Q$-Fano varieties are
rationally  connected}\endheading

Let $X$ be a log $Q$-Fano variety, i.e, if there
exists an effective $Q$-divisor $D$ such that
 the pair $(X,D)$ is Kawamata log terminal ({\it klt}) and $-(K_X+D)$ is nef
and big. By a result of Miyaoka-Mori [15], $X$ is uniruled.
 A standard conjecture ([10], [12], [13], [16]) predicts that  $X$ is actually rationally
 connected. In this paper we 
apply the theory of  weak (semi) positivity of the direct images of
 (log) relative
dualizing sheaves $f_{*}(K_{X/Y}+\Delta)$ (which has
been developed by 
 Fujita, Kawamata, Koll\'{a}r, Viehweg and others) to show 
that  a log $Q$-Fano variety is indeed
rationally  connected.

{\it Remark: } The rational connectedness of   smooth Fano varieties was established by 
Campana [1] and Koll\'{a}r-Miyaoka-Mori [12]. However their approach
 relies heavily   on the (relative) deformation theory
which seems   quite difficult to extend to the singular
case.

\proclaim{Theorem 1} Let $X$ be a log $Q$-Fano
variety. Then $X$ is rationally  connected, i.e., for
any two 
closed points  $x,y\in X$ 
there exists a rational  curve $C$ which contains $x$
and $y$. 
\endproclaim

{\it Remark: } When $\text{dim}(X)\leq 3$ and $D=0$, this was
 proved  by Koll\'{a}r-Miyaoka- Mori in [13].
 On the other hand, the result is false  for log canonical
singularities (see 2.2 in [13]).

As a corollary, we can show the following result which
 was obtained by S. Takayama [16].

\proclaim{Corollary 1} Let $X$ be a log $Q$-Fano
variety.Then $X$ is simply connected.
\endproclaim

\demo{Proof} $\pi_{1}(X)$ is finite by Theorem 1 and a
result of F. Campana [2]. On the other hand, 
$h^{i}(X,\Cal O_X)=0$ for $i\geq 1$ by
 Kawamata-Viehweg vanishing ([4], [9]). Thus we 
have $\chi(X, \Cal O_X)=1$ and hence $X$ must be
simply connected.
\enddemo

Our Theorem 1 is a consequence of the following
proposition.

\proclaim{Proposition 1} Let $X$ be a log $
Q$-Fano variety and let $f:X--->Y$ be a dominant
rational map, where $Y$ is a projective variety. Then
$Y$ is uniruled if $\text{dim}Y>0$.
\endproclaim

Let us assume  Proposition 1. Let $X'$ be a
resolution of $X$, then $X'$ is uniruled. There exists a
nontrivial maximal rationally
 connected fibration $f:X'--->Y$ ([1], [12]).
 By a result of Graber-Harris-Starr
 [5],  $Y$ is not uniruled.
 However   Proposition 1 tells us that $Y$ must be
a point and hence $X$ is rationally  connected.

The general  strategy  for proving Proposition 1 is as follows.
By a result of Miyaoka-Mori [15], it suffices to construct a covering family
of curves $C_{t}$ on $Y$ with $C_{t}\cdot K_Y<0$ for every $t$. To this end, we 
  apply the positivity theorem of the direct images of (log) relative
dualizing sheaves to $f:X--> Y$.  We  show that there  exist an ample $Q$-divisor $H$ on $Y$
and  an effective 
$Q$-divisor $D$ on $X$ such that $-K_Y=D+f^{*}H$ (modulo some exceptional
divisors). Now let $C_{t}=f(F_t)$, where $F_t$ are the general complete intersection curves
 on $X$. We have a covering family
of curves $C_{t}$ on $Y$ with $C_{t}\cdot K_Y<0$ for every $t$.

 Before we start to prove Proposition 1, let us
first give some related definitions. Also the proof
of Proposition 1 depends heavily on Kawamata's paper
[8] and Viehweg's paper [17]. 

We work over the complex number field $\Bbb C$ in  this paper.

\proclaim {Definition 1 [9]} Let $X$ be a normal
projective variety
of dimension $n$ and $K_X$  the canonical divisor on
$X$. Let $D=\sum a_{i}D_{i}$ be an effective
$Q$-divisor on $X$, where $D_{i}$  are distinct
irreducible divisors and  $ a_{i}\geq 0$. The
 pair $(X,D)$ is said to be Kawamata log terminal ({\it
klt})
(resp.  log canonical) if  $K_X+D$ is a
$Q$-Cartier divisor and if there exists a
desingularization ({\it log resolution})
$f:Z\rightarrow X$ such that the union $F$ of the
exceptional locus of $f$ and the inverse image of the
support of $D$ is a divisor with normal crossing and
$$ K_{Z}=f^{*}(K_{X}+D)+\sum_{i} e_{j}F_{j},$$ with
$e_{j}>-1$ (resp. $e_{j}\geq -1$).  $X$ is
 said to be Kawamata 
log terminal (resp. log canonical) if so is  $(X,0)$.
\endproclaim
\proclaim {Definition 2}  Let $X$ be a normal
projective variety
of dimension $n$ and $K_X$  the canonical divisor on
$X$. We say $X$ is a $Q$-Gorenstein variety if there
exists some
integer $m>0$ such that $mK_X$ is a Cartier divisor. A
$Q$-Cartier divisor $D$ is said to be {\it nef} if the
intersection number $D\cdot C\geq 0$ for any curve $C$
on $X$. $D$ is said to be {\it big} if the
Kodaira-Iitaka dimension $\kappa (D)$ attains the maximum $\text{dim}X$.
\endproclaim

The following lemma due to Raynaud [19] is quite
useful:

\proclaim {Lemma 1} Let $g:T\rightarrow W$ be
surjective morphism
of smooth varieties. Then there exists a birational
morphism of smooth
variety $\tau: W'\rightarrow W$ and a
desingularization $T'\rightarrow
T\times_{W} W'$, such that the induced morphism
$g': T'\rightarrow W'$ has the following property: 
 Let $B'$ be any
divisor of $T'$ such that  $\text{codim}(g'(B'))\geq
2$.Then $B'$ lies
in the exceptional locus of $\tau: T'\rightarrow T$.
\endproclaim

\demo{Proof of  Proposition 1} By the Stein
factorization and  desingularizations, we may assume
that $Y$ is smooth. Resolving  the indeterminacy of
$f$ and
taking a log resolution, we have a smooth projective
variety $Z$ and the surjective morphisms
$g:Z\rightarrow Y$ and $\pi:Z\rightarrow X$, $$
\CD
Z  @>\pi>> X \\
@VVgV         @. \\
Y   @.     
\endCD
$$ such that
$$ K_{Z}={\pi}^{*}(K_{X}+D)+\sum_{i} e_{i}E_{i}$$ with
$e_{i}>-1$, where $\sum E_{i}$ is a divisor with 
normal crossing.

 Since $-(K_X+D)$ is nef and big, by Kawamata
base-point free theorem, there exists an effective 
 $Q$-divisor $A$ on $X$ such that $-(K_X+D)-A$ is an ample $
Q$-divisor. Thus we can choose another ample $Q$-divisor
 $H$ on $Y$ and an effective
 $Q$-divisor $\Delta=\sum_{i} \delta_{i}E_i$ on $Z$
(with small $\delta\geq 0$ and $E_i$ are
$\pi$-exceptional if 
  $\delta_i>0$) such that
$-\pi^{*}(K_X+D+A)-\Delta-g^{*}(H)=L$ is again an ample $Q$-divisor on
$Z$. We may also assume that $\text{Supp}(\Delta+\cup_{i}
E_i+A)$ is a 
divisor with simple normal crossing
and the pair
 $(X,L+\Delta+\sum_{i} \{-e_{i}\}E_{i}+A)$ is klt. Thus we
have
$$K_{Z/Y}+\sum_{i} {\epsilon}_{i}E_{i}\sim_{Q} \sum_{i} 
m_{i}E_{i}-g^{*}(K_Y+H)$$ where $m_i={\lceil e_{i}
\rceil}$ are non-negative
integers ($\{,\}$ is the fractional part and
$\lceil, \rceil$ is the round up). 
 Also $E_i$ on the left side  contain components of
$L$,
$A$ and $\Delta$  with $0\leq \epsilon_i<1$.

We follow closely from Kawamata's paper [8].

After further blowing-ups if necessary, and  by Lemma
1, we may assume
that:
\roster
\item  There exists a  normal crossing divisor
$Q=\sum_{l} Q_{l}$  on $Y$ such that $g^{-1}(Q)\subset
\sum_{i} E_{i}$ and $g$ is smooth over $Y\setminus Q$. 
\item  If a divisor $W$ on $Z$ with
$\text{codim}(g(W))\geq 2$, then $W$ is
 $\pi$-exceptional.
\endroster

Let $D=\sum_{i} ({\epsilon}_{i}-m_{i})E_{i}=\sum_{i}
{d}_{i}E_{i}=D^{h}+D^{v}$, 
where 
\roster
\item $g: D^h\rightarrow Y$ is
surjective and smooth over 
$Y\setminus Q$ (we say $D^h$ is $g$-horizontal)
\item $g(D^v)\subset Q$ (we say $D^v$ is
  $g$-vertical [8]).
\endroster

 Notice that here besides those $E_i$  from the log
resolution
$\pi$, $D$ also contains the components of $\text{Supp}
(L+\Delta+A)$. Nevertheless ${d}_{i}<1$ for all $i$.

We have two cases:

(1):  Every $\pi$-exceptional divisor $E_i$ with
${d}_i<0$ is $g$-vertical. In this case, the natural 
homomorphism $\Cal O_Y\rightarrow g_{*}({\lceil -D
\rceil})$ is
surjective at the generic point of $Y$. 

Let
$$g^{*}Q_{l}=\sum_{j}w_{lj}E_{j}$$ and 
$$a_{j}=\frac{{d}_{j}+w_{lj}-1}{w_{lj}} \text { if } g(E_j)=Q_{l}$$ and
$$b_{l}=\{\text{max}\{a_j\}: g(E_j)=Q_{l}\}.$$

Let $N=\sum_{l} b_{l}Q_l$ and $M=-H-K_Y-N$.
Then by a result of Kawamata [8, Theorem 2], $M$ is
nef.
 
On the other hand, $g^{*}N=F+G$, where
\roster
\item 
$\text{Supp}(F)$ is $\pi$-exceptional (from those
$Q_l$ with $g^{*}Q_{l}=\sum_{j}w_{lj}E_{j}$
and each $E_j$ is $\pi$-exceptional). 
\item $G$ is effective (from those $Q_l$ with 
$g^{*}Q_{l}=\sum_{j} w_{lj}E_{j}$
 and at least one $E_j$ has coefficient $d_j\geq 0$).
\endroster

Now let $C$ be a general complete intersection curve
on $Z$ such that $C$ does not intersect with
$\text{Supp}(F)$
(e.g., pull-back a general complete intersection curve
from $X$).
Then $$g_{*}(C)\cdot (-K_Y-H)=C\cdot
g^{*}(-K_Y-H)=C\cdot g^{*}(M+N)\geq 0$$ and hence
$g_{*}(C)\cdot K_Y\leq g_{*}(C)\cdot (-H)<0$. Since
the family of the  curves $g(C)$ covers a Zariski 
open set of $Y$, by [15] $Y$ is uniruled.
 
(2):  Some $\pi$-exceptional divisor is $g$-horizontal, in
particular, $D^h$ is not zero.
 
 We have $K_{Z/Y}+D\sim_{Q} -g^{*}(K_Y+H)$. By the
stable reduction theorem and the covering trick
([8],[17]), there exists a finite morphism
$p:Y'\rightarrow Y$ such 
that $Q'=\text{Supp}(p^{*}Q)$
is a normal crossing divisor and the induced morphism
$g':Z'\rightarrow Y'$ from a desingularization 
$Z'\rightarrow Z\times_{Y}Y'$
is semistable over $Y'\setminus B$ with $\text{codim}(B)\geq
2$. Let $Z'\rightarrow Z$ be the induced morphism.

$$
\CD
Z'    @>q>> Z  @>\pi>> X \\
@Vg'VV   @VgVV   @. \\
Y'    @>p>>    Y @.
\endCD
$$

 We can write  
$$K_{Z'/Y'}+D'\sim_{Q}g'{*}p^{*}(-K_Y-H),
\text{where } D'=\sum_{j'} d'_{j'}E'_{j'}.$$

The coefficients $d'_{j'}$ can be calculated as follows [8]:

\roster
\item If $E'_{j'}$ is $g'$-horizontal and
$q(E'_{j'})=E_j$, then $d'_{j'}=d_j$.
\item If $E'_{j'}$ is $g'$-vertical with $q(E'_{j'})=E_j$ and 
$g'(E'_{j'})=Q'_{l'}$, then
$d'_{j'}=e_{j'}(d_{j}+w_{lj}-1)$, where $e_{l'}$ is
the ramification index of
$q$ at the generic point of $E'_{j}\rightarrow E_j$.
\item We are not concerned with those 
$E'_{j}$ such that $g'(E'_{j})\subset B$. 
\endroster

{\it Note}: In Kawamata's paper [8], he replaced $D$
by $D-g^{*}N$. Thus
all the coefficients there $\text{are}\leq 0$.
However, here we do not make such replacement
and in our case some coefficients
$d_{i}={\epsilon}_{i}>0$

 Therefore, by the standard trick (keep the fractional part on the
 left side and the integral part on the right side, also blow-up $Z'$
if necessary). We have 

$$K_{Z'/Y'}+\sum_{j'} {{\epsilon}'}_{j'}E'_{j'}\sim_{Q}
\sum_{k'} n_{k'}E'_{k'}-{g'}^{*}p^{*}(K_Y+H)-V+G$$ such that
\roster
\item  $\sum_{k'} n_{k'}E'_{k'}$ is Cartier, $\text{Supp}(\sum_{k'} E'_{k'})$ 
 is  $q\circ\pi$-exceptional and $g'(\sum_{k'} n_{k'}E'_{k'})$ is not contained
in $B$.
\item $V$ is an effective Cartier divisor which is  $g'$-vertical (from those $g$-vertical $E_i$ with
$d_i=\epsilon_i\geq 0$ and $d'_{j'}=e_{j'}(d_{j}+w_{lj}-1)\geq 1$ for some $l$.)
\item $G$ is also Cartier and $q\circ\pi$-exceptional (from those
$E'_{j}$ with $g'(E'_{j})\subset B$).
\item $\sum_{j'} d_{j'}E'_{j'}$ remains on the left side, where $E'_{j'}$
are $g'$-horizontal with $0<d_{j'}=d_j={\epsilon}_{j}<1$
\item If $E'_{k'}$ is $g'$-horizontal, then
$n_{k'}\geq 0$. 
\item  $(Z',\sum_{j'} {{\epsilon}'}_{j'}E'_{j'})$ is klt.
\endroster 

There also exists a cyclic cover [9] $p': Y''\rightarrow
Y'$ such that $Y''$ is smooth and
${p'}^{*}p^{*}(H)=2H'$, where $H'$
is an ample Cartier divisor. Since $H$ is  an ample $Q$-divisor, we can
choose  the covering 
 $p'$ in such way that the ramification locus $R_{p'}$
of $p'$ intersects $\text{Supp}Q$ and $B$
 transversely. Let 
$g'':Z''\rightarrow Y''$ be the induced morphism from
a desingularization 
$Z''\rightarrow Z'\times_{Y'}Y''$.

$$
\CD
\
Z''  @>q'>>Z'    @>q>> Z  @>\pi>> X \\
@Vg''VV  @Vg'VV   @VgVV   @. \\
Y''  @>p'>>Y'    @>p>> Y   @.
\endCD
$$

 Since $g'$ is semistale over $Y'\setminus B$, we have
$q'{*}K_{Z'/Y'}=K_{Z''/Y''}$ over $Y''\setminus {p'}^{-1}(B)$.
  Thus  again we can
 write $$K_{Z''/Y''}+\sum_{j''}
{{\epsilon}''}_{j''}E''_{j''}\sim_{Q}
\sum_{k''} n_{k''}E''_{k''}-{g''}^{*}{p'}^{*}p^{*}K_Y-2{g''}^{*}H'-V'+G'$$
where 
\roster
\item $\sum_{k''} n_{k''}E''_{k''}$ is Cartier and
 $\text{Supp}(\sum_{k''}  E''_{k''})$ is ${q'\circ
q}\circ \pi$-exceptional.
\item $V'$ is an effective Cartier divisor which is  $g''$-vertical.
\item $G'$ is also Cartier and ${q'\circ q}\circ \pi$-exceptional
(since $\text{codim}g''(G')\geq 2$).
\item  If $E''_{k''}$ is $g''$-horizontal, then
$n_{k''}\geq 0$. 
\item $(Z'',\sum_{j''} {{\epsilon}''}_{j''}E'_{j''})$ is klt
and $\sum_{j''} {{\epsilon}''}_{j''}E''_{j''}$ is $Q$-linearly equivalent to 
a Cartier divisor.
\endroster

Since all the
$g''$-horizontal divisors $E''_{k''}$ have non-negative
coefficients,
 $${g''}_{*}(\sum_{k''}
n_{k''}E''_{k''}-{g''}^{*}{p'}^{*}p^{*}K_Y-2{g''}^{*}H'-V'+G')$$
is not a zero sheaf. Let $\omega =\sum_{k''}
n_{k''}E''_{k''}-{g''}^{*}{p'}^{*}p^{*}K_Y-V'+G'$.

 Applying  the results of Koll\'{a}r [11], Viehweg [17, Lemma 5.1] and Kawamata [7, Theorem 1.2], we may
assume that $${g''}_{*}(\omega -2{g''}^{*}H')={g''}_{*}(\omega)\otimes \Cal
O_{Y''}(-2H')$$ is torsion
free and weakly positive over $Y''$.

{\it Note: } In [7], Kawamata proved that in fact (after blow-up $Y''$ further) we may assume that
${g''}_{*}(\omega )$ is locally free and semipositive. However the weak positivity
is sufficient in our case.

 By the weak positivity of ${g''}_{*}(\omega)\otimes \Cal
O_{Y''}(-2H')$, we have 
$$\Cal
{\hat {\Cal S}}^{n}{g''}_{*}(\omega)\otimes \Cal
O_{Y''}(-2nH'+nH')$$is 
 generically generated by its global sections over $Y''$ for some $n>0$, where ${\hat {\Cal
S}^{n}}$ denotes 
the reflexive hull of ${\Cal S}^{n}$.

We have a natural homomorphism $$S^{n}{g''}^{*}{g''}_{*}(\omega)\otimes 
{g''}^{*}\Cal
O_{Y''}(-nH')\rightarrow \omega^{n}\otimes {g''}^{*}\Cal O_{Y''}(-nH')$$

By the torsion freeness of ${g''}_{*}(\omega)$, there exists
an open set $U\subset Y''$ with $\text{codim}(Y''\setminus U)\geq 2$ such that
$$\hat {\Cal S}^{n}{g''}_{*}(\omega)\otimes \Cal
O_{Y''}(-nH')|_{U}={\Cal S}^{n}{g''}_{*}(\omega)\otimes \Cal
O_{Y''}(-nH')|_{U}$$ and hence $${g''}^{*}{\hat {\Cal S}}^{n}{g''}_{*}
(\omega)\otimes {g''}^{*}\Cal
O_{Y''}(-nH')|_{W}=S^{n}{g''}^{*}{g''}_{*}(\omega)\otimes {g''}^{*}\Cal
O_{Y''}(-nH')|_{W}$$ where $W={g''}^{-1}(U)$. If $B'=Z''\setminus W$,
 then $B'$ is $g''$-exceptional.
Since ${g''}^{*}\hat {\Cal S}^{n}{g''}_{*}(\omega)\otimes {g''}^{*}\Cal
O_{Y''}(-nH')$ is also generically generated by its global sections over $Z''$,
there is a non trivial morphism $$\bigoplus \Cal O_{Z''}|_{W}\rightarrow
 \omega^{n}\otimes {g''}^{*}\Cal O_{Y''}(-nH')|_{W}$$ i.e., $\omega^{n}\otimes 
{g''}^{*}\Cal O_{Y''}(-nH')$ admits a meromorphic section which has poles only along
$B'$. Thus we may choose some large integer $k$ such that $\omega^{n}\otimes 
{g''}^{*}\Cal O_{Y''}(-nH')+kB'$ has a holomorphic section, i.e.,   $$n(\sum_{k''}
n_{k''}E''_{k''}-{g''}^{*}{p'}^{*}p^{*}K_Y-V'+G')-n{g''}^{*}H'+kB'$$ is effective. 

Again as before, we can choose a  family of general
complete intersection curves $C$ on $Z''$ such that
$C$ does not intersect
with the exceptional locus of $Z''\rightarrow X$
(such as ${E''}_{k''}$, $B'$ and $G'$). Thus
${g''}_{*}(C)\cdot {p'}^{*}p^{*}(K_Y)\leq {g''}_{*}(C)
\cdot (-H')<0$ and hence $Y$ is uniruled [15]. q.e.d.
\enddemo

\heading{\bf $\S$2.The behavior of the canonical
bundles under projective 
morphisms}\endheading

Let $X$ and $Y$ be two projective varieties and
$f:X\rightarrow Y$ be a surjective morphism. Assume
that the Kodaira dimension $\kappa (X)\leq 0$. In
general, it is almost impossible to predict the
Kodaira
dimension of $Y$. The following example shows
that even when  $\text{dimY}=1$, we have  no control
of the genus of $Y$:

\proclaim{Example} ([6], [13]):  Let $C$ be a smooth curve of
arbitrary genus
$g$ and $A$ be an ample line bundle on $C$ such that
$deg A>2deg K_{C}$. Let $S=\text{Proj}{_C}({\Cal
O}(A)\oplus {\Cal 
O}_{C})$
be the projective space bundle associated to the
vector bundle $\Cal O(A)\oplus {\Cal O}_{C}$. Then 
$K_{S}={\pi}^*(K_{C}+A)-2L$, where $L$ is the 
tautological bundle and $\pi:S\rightarrow C$ is the
projection. An easy computation show that $-K_S$ is
big (i.e, $h^{0}(S,-mK_S)\approx c\cdot m^{2}$ for
$m\gg 0$ and some 
$c>0$).  In
particular, the Kodaira dimension $\kappa
(S)=-\infty$. However, the genus of $C$ could be 
large.\endproclaim

  On the other hand, it is not hard to find   the 
following facts:
\roster
\item Let $D=-K_S=2(L-{\pi}^{*}A)+{\pi}^{*}(A-K_C)$,
then $D$ is effective and the pair $(S,D)$ is not
 log canonical [9].
\item $-K_S$ is not nef, i.e, there exists some curve
$B$ (e.g., choose $B$ to be the section corresponding to
${\Cal O}(A)\oplus{\Cal O}_{C}\rightarrow {\Cal
O}_{C}\rightarrow 0$)
 on $S$ such that $(-K_{S})\cdot B<0$.
\item For any integer $m>0$, the linear system $-mK_S$
contains some fixed component (e.g.,
$m(L-{\pi}^{*}A)$) which
dominates $C$.\endroster

In view of the above example, we
 give  a few  sufficient conditions  which 
guarantees  nice behavior of the   Kodaira 
dimension (and the canonical bundle).

\proclaim{Theorem 2} Let $f: X\rightarrow Y$ be a
surjective
morphism. Assume that $D\equiv -K_X$ is an
 effective $Q$-divisor and the pair
$(X,D)$ is  log canonical. Moreover assume
that $Y$ is normal and $Q$-Gorenstein. Then either $Y$
is uniruled or   $K_Y$ is numerically trivial. (In
particular, $\kappa (Y)\leq 0$.)
\endproclaim
\proclaim{Corollary 2}  Let $f: X\rightarrow Y$ be a
surjective
morphism.  Assume that $X$ is log canonical
 and $K_X$ is numerically trivial 
(e.g., a Calabi-Yau manifold). Moreover assume that
$Y$ is normal and
$Q$-Gorenstein. Then either $Y$ is uniruled or  $K_Y$
is numerically trivial.
\endproclaim

\proclaim{Theorem 3} Let $f: X\rightarrow Y$ be a
surjective
morphism. Assume that  $-K_X$ is nef and $X$ is
smooth. Moreover 
assume that $Y$ is normal and
$Q$-Gorenstein. Then either $Y$ is uniruled or  $K_Y$
is numerically trivial.
\endproclaim

{\it Remark}: Theorem 3  was  proved in
[18] (in particular, the result solved a conjecture
proposed by Demailly, Peternell and Schneider [3]).
However, the proof given there was incomplete (as
pointed  out  to me by  Y. Kawamata. I wish to thank him). 
The problem lies in
Proposition 1 in [18], the  point is that the
nefness in general is
not preserved under the deformations (mod $p$ reductions in our case). We shall present
a new proof of  this proposition (see Proposition 2).

\proclaim{Theorem 4} Let $f: X\rightarrow Y$ be a
surjective
morphism. Assume that there exists some integer $m>0$
such that $-mK_X$ is effective and has no fixed locus
which dominates
$Y$. Moreover assume that $X$ is  log canonical and
$Y$ is
normal 
and
$Q$-Gorenstein. Then either $Y$ is uniruled or  $K_Y$
is numerically trivial.
\endproclaim

As an immediate consequence,  we can show the
following result about the Albanese maps.

\proclaim{Corollary 3} Let $X$ be a smooth projective
variety. Then the
Albanese map $\text{Alb}_{X}:X\rightarrow
\text{Alb}(X)$  of  $X$ is
surjective and has connected fibers if 
$X$ satisfies one of the following conditions:
\roster
\item $D\equiv -K_X$ is an effective $Q$-divisor and
the
pair $(X,D)$ is log canonical (``$\equiv$''
means {\it numerically equivalent}). 
\item $-K_X$ is nef .
\item There exists some integer $m>0$ such that
$-mK_X$ is effective and 
has no fixed component which dominates
$\text{Alb}(X)$. \endroster 
\endproclaim

The main ingredients of the proofs are the Minimal
Model Program (in particular, a vanishing theorem of
Esnault-Viehweg, Kawamata and Koll\'{a}r plays an
essential role), and the deformation theory. It is 
interesting to notice that
the proofs of Theorem 2 and Theorem 3 are completely 
 different in nature. 

The following vanishing theorem of  Esnault-Viehweg,
Kawamata and 
Koll\'{a}r is important to us:

\proclaim{Vanishing Theorem [4], [22]} Let $f:X\rightarrow
Y$ be a
surjective morphism from a smooth projective variety
$X$ to a normal variety $Y$. Let $L$ be a line bundle
on $X$ such that $L\equiv f^{*}M+D$, where $M$ is a
$Q$-divisor on $Y$
and $(X,D)$ is Kawamata log terminal. Let $C$ be a
reduced divisor without
common component with $D$ and $D+C$ is a normal
crossing divisor. Then\roster
\item $f_{*}(K_X+L+C)$ is torsion free [20].
\item Assume in addition that $M$ is nef and big. Then
$H^{i}(Y, {\Cal R}^{j}f_{*}(K_X+L+C))=0$ for $i>0$
and
$j\geq 0$.\endroster
\endproclaim

{\it Remark}: The $C=0$ case was done by 
Esnault-Viehweg, Kawamata 
and Koll\'{a}r [11]. The generalized version given here
essentially was proved
by Esnault-Viehweg in [4] and by Fujino in [22]. I thank Professor
Viehweg for informing on the matter.
 C. Hacon pointed out an inaccuracy  and informed me 
 of the reference [22]. I would like to thank him. Below, we give an outline of the proof which was
provided to me by E. Viehweg.

\demo{\it Sketch of the proof} By [4, 5.1 and 5.12],
we have an injective
morphism $$0\rightarrow H^{j}(X,K_X+L+C)\rightarrow
H^{j}(X,K_X+L+C+B)$$ for any $j$, where $B=f^{*}(F)$
for some divisor $F$ on $Y$. If we choose $F$ to be a
very ample divisor, we have
the exact sequence:$$0\rightarrow {\Cal
R}^{j}f_{*}(K_X+L+C)\rightarrow  {\Cal
R}^{j}{f_*}(K_X+L+C+B)
\rightarrow {\Cal R}^{j}f_{*}(K_B+L+C)\rightarrow 0.$$
By
induction on $\text{dim} Y$ and the
 Leray-spectral sequence associated with  ${\Cal
R}^{j}f$, we can prove the result (see [4]
for details). q.e.d.\enddemo

\demo{\it Proof of Theorem 2}: Let  $g: Z\rightarrow
X$ be a
log resolution and let $\pi=f\circ g$. Then $$K_Z=g^{*}(K_X+D)+\sum
a_{i}E_{i}, \text{ where each } a_{i}\geq -1.$$ We can
rewrite $\sum a_{i}E_i=\sum b_{j}E_j+\sum c_{k}E_k +\sum d_{l}E_l$
where $b_{j}\geq 0$, $0> c_{k}>-1$ and
$d_{l}=-1$. 

Let $C$ be a general complete intersection curve on
$Y$ and
 $W={\pi}^{-1}(C)$. We have 
$$K_Z=K_{Z/Y}+{\pi}^*K_Y \text { and }
K_{Z/Y}|_W=K_{W/C}$$  

Thus $$K_W+{\pi}^{*}(K_{Y}|_C)+\sum
-c_{k}{E_{k}}|_W+\sum  \{-b_j\}E_j|_W+\sum
{E_{l}}|_W\equiv {\pi}^{*}K_C+\sum {\lceil
b_j\rceil}E_j|_W$$ 

Let us assume that  $K_{Y}\cdot C>0$. Since  $(K_W, \sum
-c_{k}{E_{k}}|_W+\sum \{-b_j\}E_j)$ is Kawamata log
terminal and $\sum {\lceil b_j\rceil}E_j|_W$ is
exceptional, the Vanishing Theorem yields 
$$H^{1}(C,{\pi}_{*}({\pi}^{*}K_C+\sum
{\lceil b_j \rceil}E_j|_W))=H^{0}(C,{\Cal O}_C)=0,$$
a contradiction. So we must have $K_{Y}\cdot C\leq 0$.
If $K_{Y}\cdot C<0$, $Y$ is uniruled by [15]. If
$K_Y\cdot C=0$, $K_Y$ is numerically trivial by Hodge
index theorem. q.e.d.
\enddemo

\demo{\it Proof of Theorem 3} Let us first establish
the following 
proposition (Proposition 1 in [18]). \enddemo

\proclaim{Proposition 2} Let $\pi:X\rightarrow Y$  be
a surjective 
morphism  between  smooth projective varieties over 
$\Bbb C$. 
Then  for any ample divisor $A$ on $Z$, 
$-K_{X/Y}-\delta \pi^{*}A$ is 
not 
nef for any $\delta>0$.
\endproclaim

\demo{\it Proof of Proposition 2} We shall  give a new proof of this
proposition by modifying the arguments we 
  used before [18]. Again, the main idea and method comes
from [14].

Let $C\subset X$ be a general smooth curve of genus
$g(C)$ such that
 $C\nsubseteq \text{Sing}(\pi)$. Let  $p\in C$ 
 be a general point and $B=\lbrace p \rbrace$ be the
base scheme. 
Denoting by $\nu: C\rightarrow  X$ the embedding of
$C$ to $X$.
 Let $D_{Y}(\nu, B)$ be  the Hilbert scheme
representing the functor of 
the 
relative deformation 
over $Y$ of $\nu$. Then by [14], we have
$$\text{dim}_{\nu}D_{Y}(\nu, B)
\geq -\nu_{*}(C)\cdot K_{X/ Y}-g(C)\cdot \text{dim}
X$$ Suppose now 
that 
$-K_{X/Y}-\delta \pi^{*}A$ is nef for some $\delta>0$.
Let $H$ be an 
ample 
divisor on $X$ and $\epsilon>0$ be  a small number, we
may assume that
 $$\nu_{*}(C)\cdot (\delta \pi^{*}A-\epsilon H)>0.$$
Then
$$\text{dim}_{\nu}D_{Y}(\nu, B)
\geq \nu_{*}(C)\cdot (\delta \pi^{*}A-\epsilon
H)-g(C)\cdot \text{dim} 
X$$
 Since $-K_{X/Y}-\delta \pi^{*}A+\epsilon H$ is an
ample divisor on $X$ 
and the 
ampleness is indeed an open  property in nature. By
the method  of 
modulo $p$  reductions [14], after composing $\nu$ with
suitable 
Frobenius morphism if necessary, we can assume that
there exists 
another
 morphism [14]
${\nu}': C\rightarrow X$ such that
\roster 
\item $\text{deg}{\nu}'_{*}(C)<\text{deg}\nu_{*}(C)$, 
where $\text{deg}\nu_{*}(C)=\nu_{*}(C)\cdot H$.
\item $\pi\circ \nu'=\pi\circ \nu$.
\endroster
 However, we have $$\nu'_{*}(C)\cdot (\delta
\pi^{*}A-\epsilon H)>
\nu_{*}(C)\cdot (\delta \pi^{*}A-\epsilon H)$$ by (1)
and (2).
This  guarantees the existence of
 a non-trivial relative deformation of $\nu'$. Since  
$\text{deg}\nu'_{*}(C)<\text{deg}\nu_{*}(C)$, this
process
 must terminate, which is absurd. q.e.d.

{\it Proof of Theorem 3  continued}: We keep the same
notations as in 
Theorem 2.
 Let $-K_{W/C}=-K_{X}|_W+f^{*}(K_{Y}|_C)$,
where $C$ is a general complete intersection curve on
$Y$ and 
$W=f^{-1}(C)$.
Applying  Proposition 2, we deduce that $K_{Y}\cdot
C\leq 0$ and 
 we 
are done. q.e.d.\enddemo

\demo{\it Proof of Theorem 4} Replacing $X$ by a suitable
resolution if
necessary, we may assume that $-pK_X=L+N$ for some positive integer $p$, where $|L|$
is base-point free and $N$ is the fixed part.
Multiplying both sides by some large integer $m$, we
can write $-K_X=\epsilon L_{m}+N_{m}$ as $Q$-divisors,
where 
$\epsilon>0$
is a small rational number. We may again assume that
the linear system $L_m$ is base-point free and $N_m$
is 
the fixed locus. The point is that $(X,\epsilon L_m)$
is
 log canonical. If the fixed locus does not
dominate $Y$, we can choose a general complete
intersection curve $C$ on $Y$ such that $C$  only
intersects
 $f(\text{Supp}(N_m))$ at some isolated points. Using the same
notations as before,
we have  $$K_W+f^{*}(K_{Y}|_C) +\epsilon
L_m|_W\equiv f^{*}K_C+E-N_m|_W$$
where $E$ is exceptional. If $C\cdot K_Y>0$, then by the 
Vanishing Theorem   $$H^{1}(C, f_{*}(K_W+f^{*}(K_{Y}|_C) +\epsilon
L_m|_W+\{N_m\}))=H^{0}(C,-f_{*}(-{\lfloor N_m|_W \rfloor}))=0$$ 
Since $\text{Supp}(N_m)$ is contained in
some fibers of $f$, we reach a contradiction.  
  The remaining
arguments are exactly
the same as in the proof of Theorem 2. q.e.d.\enddemo

\demo{\it Proof of Corollary 3} Let  $X\overset
f\to\rightarrow Y\overset 
g\to\rightarrow \text{Alb}(X)$ be the
 Stein factorization of $\text{Alb}_X$. Then by
Theorem 2-4, we 
conclude that  
$\kappa(Y)=0$. We may assume that $Y$ is smooth, otherwise we can
take a 
desingularization 
$Y'$ of $Y$.  This will not affect our choices for the
general curve 
$C$ (since 
$Y$ is smooth 
in 
codimension 1). Therefore  $\text{Alb}(X)$ must be  an
 abelian variety 
and 
hence 
$\text{Alb}_X$ is surjective and has connected fibers
(see [18] for 
details). 
q.e.d.
\enddemo

{\it Remark}: The notion of special varieties was introduced and
 studied  by F. Campana in [20]. He also conjectured that compact
K\"{a}hler manifolds with $-K_X$ nef are special. S. Lu [21] proved the
 conjecture
for projective varieties. In particular, if $X$ is a
projective variety with $-K_X$ nef and if there is a surjective map $X\to Y$. 
Then $\kappa (Y)\leq 0$. Our  focus  however, is on the uniruledness of $Y$.

 {\bf Acknowledgments}:  Part of  the results (Theorem 2-4) were
obtained while I stayed at the University of Tokyo, Japan and
 the University of Nancy, France. I would
like to thank Professor Y. Kawamata, Professor K.
Oguiso and  Professor  F. Campana for their  hospitality. I  thank both institutions for
the  financial support.  I am  grateful  to 
Professor Y. Kawamata, Professor J. Koll\'{a}r and Professor E. Viehweg  for encouragement, helpful advice,
 and  valuable suggestions. I am  grateful to Qihong Xie who has carefully read an earlier 
version of the paper.  His useful comments and suggestions helped to improve
 the paper. Finally, I  would like to thank the referee for valuable suggestions.

\enddocument